\documentclass[12pt]{amsart}
\usepackage{latexsym}           
\usepackage{amssymb}
\usepackage{amsxtra}
\usepackage{lscape}
\usepackage{epsfig}             
\usepackage{amsfonts}

\newcommand\beq{\begin{equation}}
\newcommand\eeq{\end{equation}}

\newcommand{\IP}{\mathbb{P}}                                     
\newcommand{\IQ}{\mathbb{Q}}                           
\newcommand{\IR}{\mathbb{R}}                           
\newcommand{\IC}{\mathbb{C}}

\newcommand{\cC}{\mathcal{C}}
\newcommand{\bcC}{\boldsymbol{\mathcal{C}}}
\newcommand{\bcO}{\boldsymbol{\mathcal{O}}}

\newcommand{\F}{\mathcal{F}}
\newcommand{\cF}{\mathcal{F}}

\newcommand{\cO}{\mathcal{O}}

\newcommand{\g}{       \mathfrak{g}     }

\newcommand{\PVI}{$\text{P}_{\text{VI}}$}   

\newcommand{\lsl}{\mathfrak{sl}} 

\newcommand{\fus}{\circledast}

\newcommand{\bM}{{\bf M}}

\newcommand{\pf}{\begin{bpf}}

\newcommand{\pfms}{\begin{bpfms}}
\newcommand{\epf}{\end{bpf}\hfill$\square$\\}           
\newcommand{\epfms}{\end{bpfms}\hfill$\square$\\}               

\newcommand{\idea}{\begin{bidea}}

\newcommand{\eidea}{\end{bidea}\hfill$\square$\\}           

\newcommand{\sk}{\begin{bsk}}    

\newcommand{\esk}{\end{bsk}\hfill$\square$\\}           
\newcommand{\sketch}{\begin{bsketch}}

\newcommand{\esketch}{\end{bsketch}\hfill$\square$\\}

\newcommand{\wt}{\widetilde}

\newcommand{\al}{\alpha}

\newcommand{\be}{\beta}
\newcommand{\ga}{\gamma}
\newcommand{\de}{\delta}
\newcommand{\Ga}{\Gamma}
\newcommand{\si}{\sigma}
\newcommand{\om}{\omega}

\newcommand{\Sym}{\text{\rm Sym}}

\newcommand{\tr}{\text{\rm Tr}}

\newcommand{\Hom}{\text{\rm Hom}}

\newcommand{\SL}{\text{\rm SL}}
\newcommand{\PSL}{\text{\rm PSL}}
\newcommand{\GL}{\text{\rm GL}}

\newcommand{\SO}{\text{\rm SO}}

\newcommand{\PSU}{\text{\rm PSU}}

\newcommand{\diag}{{\text{\rm diag}}}

\newcommand{\spq}{/\!\!/}

\def\mapright#1{\smash{
        \mathop{\longrightarrow}\limits^{#1}}}

\newenvironment{eqn}
        {\begin{equation}}
        {\end{equation}}

\theoremstyle{plain}
\newtheorem {hypo}{\bf\hspace{-\parindent}Hypothesis}

\newtheorem {prop}[hypo]{Proposition}

\theoremstyle{definition}
\newtheorem {defn}[hypo]{Definition}

\theoremstyle{remark}
\newtheorem {rmk}[hypo]{Remark}

\begin{document}


\title[Some explicit solutions to the Riemann--Hilbert problem]
{Some explicit solutions to the\\ Riemann--Hilbert problem}
\author[P. P. Boalch]{Philip Boalch}
\address{\'Ecole Normale Sup\'erieure\\
45 rue d'Ulm\\
75005 Paris\\
France} 
\email{boalch@dma.ens.fr}

\subjclass[2000]{34M50, 33E17, 34M55}
\keywords{Riemann--Hilbert problem, Painlev\'e equations, algebraic solutions,
Heun equations,  tetrahedral group, octahedral group,
triangle groups, Belyi maps}

\begin{abstract}
Explicit solutions to the Riemann--Hilbert problem will be found 
realising some irreducible non-rigid local systems. The relation to
isomonodromy and the sixth Painlev\'e equation will be described.

\end{abstract}

\dedicatory{
Dedicated to Andrey Bolibruch}


\maketitle


\renewcommand{\baselinestretch}{1.02}            
\normalsize

\begin{section}{Introduction}
Unfortunately, to say 
that a particular Riemann--Hilbert problem is ``solvable'', one usually
means that  {\em there exists a solution} rather than that one is actually
able to solve the problem explicitly.

In this article we will confront the problem of explicitly 
solving the Riemann--Hilbert problem directly, for irreducible representations
(so we already know the problem is ``solvable'').
We will describe how one soon becomes embroiled in isomonodromic
deformation equations, from which it is easy to see the difficulty:
in the simplest non-trivial case the isomonodromy equations reduce to the
sixth Painlev\'e equation \PVI\ and one knows that generic solutions of \PVI\
cannot be written explicitly in terms of classical special functions. 

However there are some explicit solutions to \PVI\ and our aim will be to
write down some new solutions controlling isomonodromic deformations of 
non-rigid rank
two Fuchsian systems on the four-punctured sphere.
The cases we will study here will have 
monodromy group equal to
either the binary tetrahedral or octahedral group 
(the icosahedral case having been studied in \cite{icosa}),
or to one of the
triangle groups $\Delta_{237}$ or  $\Delta_{238}$.

Previously a tetrahedral and an octahedral solution of \PVI\  have 
been constructed by Hitchin
\cite{Hit-Octa} and (up to equivalence) independently 
by Dubrovin \cite{Dub95long}.
Moreover with hindsight we see there are three other such solutions in
the work of Andreev and Kitaev \cite{And-Kit-CMP, Kitaev-dessins}.
Here we will classify all such solutions and find an explicit solution in each 
of the new cases that appear. 

Amongst the solutions which look to be new 
(i.e. to the best of the author's knowledge have not previously appeared)
there are five octahedral
solutions including one of genus one, and two $18$ branch
genus one solutions with monodromy group $\Delta_{237}$.
The largest octahedral solution has
sixteen branches which is (currently) the largest known genus zero solution
(those with more branches in \cite{icosa} having higher genus) 
and we
will show it is equivalent to a solution with monodromy group $\Delta_{238}$.

The results of sections \ref{sn: tet} and 
\ref{sn: oct} will be of particular interest
to people interested in constructing linear differential equations with
algebraic solutions (cf. e.g. 
\cite{Katz-Algsols, Bald-Dwork-algsols, vdP-Ulmer-defg, B-vdW}).
Indeed tables $1$ and $3$ may be interpreted as the analogue
for rank two Fuchsian systems with four poles on $\IP^1$, of 
the tetrahedral and octahedral parts of
Schwarz's famous list \cite{Schwarz} 
of hypergeometric equations with algebraic bases of solutions. 

\end{section}

\begin{section}{From Riemann--Hilbert to Painlev\'e}

Consider a logarithmic connection $\nabla$ on the trivial rank $n$ 
complex vector bundle
over the Riemann sphere with singularities at points 
$a_1,\ldots,a_m$. Choosing a coordinate $z$ on the sphere (in which
$a_m=\infty$ say), this amounts to giving the Fuchsian system of
differential equations $\nabla_{d/dz}$
which will have the form: 
\begin{equation} \label{eqn: 2x2 linear}
\frac{d}{dz}- A(z); \qquad A(z)=\sum_{i=1}^{m-1}\frac{A_i}{z-a_i}
\end{equation}
for complex $n\times n$ matrices $A_i$.
The original Riemann--Hilbert map is the map which takes such a Fuchsian system
to its monodromy data:
restricting $\nabla$ to the punctured sphere 
$$\IP^*:=\IP^1\setminus\{a_1,\ldots,a_m\}$$
yields a nonsingular holomorphic connection and taking its monodromy yields a
representation
$$\rho\in\Hom(\pi_1(\IP^*),G)$$
where $G=\GL_n(\IC)$.
The Riemann--Hilbert problem is the following: 
given $a_1,\ldots,a_m$ and $\rho$ can we find 
such a connection $\nabla$ with monodromy equal to $\rho$? 

Upon choosing simple loops $\ga_i$ in $\IP^*$ around $a_i$ generating 
$\pi_1(\IP^*)$ and such that $\ga_m\circ\cdots\circ\ga_1$ is contractible
one sees that for each $m$-tuple of points ${\bf a}=(a_1,\ldots,a_m)$ 
the Riemann--Hilbert map
amounts to a map between the following spaces:
\begin{equation}\label{eq: rhmap}
 \bigl\{ (A_1,\ldots,A_m) \ \bigl\vert \ \text{$\sum A_i= 0$} \bigr\}
\ \mapright{\text{RH}_{\bf a}}\ 
\bigl\{ (M_1,\ldots,M_m) \ \bigl\vert \ M_m\cdots M_1 = 1 \bigr\}
\end{equation}
where $M_i=\rho(\ga_i)\in G$.
The Riemann--Hilbert problem then becomes:
given a point ${\bf M}=(M_1,\ldots,M_m)$ 
on the RHS of \eqref{eq: rhmap}, are there matrices 
${\bf A}=(A_1,\ldots A_m)$ with $\sum A_i=0$  
on the LHS such that $\text{RH}_{\bf a}({\bf A})={\bf M}$?

\begin{rmk}
So far we have ignored the questions of choosing a basepoint for 
$\pi_1(\IP^*)$ and the choice of basis of the fibre at the basepoint. 
However it is immediate that if we have a solution 
$\text{RH}_{\bf a}({\bf A})={\bf  M}$ (defined with respect to some choice of
basepoint/basis) then conjugating the matrices $A_i$ by
some constant matrix $g\in G$ corresponds to conjugating the monodromy
matrices $M_i$ as well. Thus the Riemann--Hilbert problem is independent of
the 
choice of basepoint/basis since these just move to conjugate representations. 
\end{rmk}

Some fundamental work on the Riemann--Hilbert problem was done by Schlesinger 
\cite{Schles}. 
He considered the question of constructing new Riemann--Hilbert solutions from
a given solution $\text{RH}_{\bf a}({\bf A})={\bf  M}$, in two ways:

\noindent $1)$ 
   Schlesinger examined the fibres of the Riemann--Hilbert map and
   defined ``Schlesinger transformations'', which move ${\bf A}$ 
   within the fibres (cf. also \cite{JM81}). Roughly
   speaking generic fibres are discrete and correspond
   to certain integer shifts
   in the eigenvalues of the matrices $A_i$; geometrically these Schlesinger
   transformations amount to rational gauge transformations with singularities
   at the poles of the Fuchsian system.

\noindent $2)$ 
Schlesinger also found how the matrices ${\bf A}$ can be varied as one moves 
the pole positions ${\bf a}$ in order to realise the same monodromy data 
${\bf M}$. (Locally---for small deformations of ${\bf a}$---this 
makes sense as one can use the same loops generating $\pi_1(\IP^*)$;
globally one should drag the loops around with the points ${\bf a}$, so
on returning ${\bf a}$ to their initial configuration 
$\rho$ may have changed by the
action of the mapping class group of the $m$-pointed sphere.)
He discovered that  if the matrices $A_i$ satisfy the
following non-linear differential equations, now known as the Schlesinger
equations, then locally the monodromy data is preserved (up to overall
conjugation): 

\begin{equation} \label{eqn: 2x2 schles}
\frac{\partial A_i}{\partial a_j}=\frac{[A_i,A_j]}{a_i-a_j}\qquad 
\text{if } i\ne j, \text{ and}\qquad
\frac{\partial A_i}{\partial a_i}=-\sum_{j\ne
i}\frac{[A_i,A_j]}{a_i-a_j}.
\end{equation}

In the generic case such an ``isomonodromic deformation'' necessarily satisfies
these equations (up to conjugation).
This gives a hint at the difficulty of the Riemann--Hilbert problem: even if
one knows a solution for some configuration of pole positions, one must
integrate some nonlinear differential equations to obtain solutions for a
deformed configuration.

This also gives a hint at how one might find some interesting solutions to the
Riemann--Hilbert problem. Namely since one can move the pole positions one may
consider degenerations into systems with fewer poles (for which the problem
should be easier). Using solutions to these degenerate Riemann--Hilbert 
problems one can get 
asymptotics for the original solution to the Schlesinger
equations and in good circumstances this enables  computation of the solution. 
This is in effect what we will do below (using the analysis of the
degenerations in \cite{SMJ} part II and \cite{Jimbo82}).

Suppose we fix an irreducible representation 
$\rho\in\Hom(\pi_1(\IP^*),G)$. 
Let $\cC_i\subset G$ be the conjugacy class containing $M_i=\rho(\ga_i)$ 
which we will
suppose for simplicity is regular semisimple, 
although this is not strictly necessary. 
(We are thus considering ``generic'' representations.)
 
Since $\rho$ is irreducible we know \cite{AnoBol94} 
there exists some Riemann--Hilbert 
solution $\text{RH}_{\bf a}({\bf A})={\bf M}$.
Let $\cO_i\subset \g$ be the adjoint orbit of $A_i$ 
(in the Lie algebra of $n\times n$ complex matrices).
By genericity we know $\exp(2\pi\sqrt{-1}\cO_i)=\cC_i$.
Indeed if in the Riemann--Hilbert map we restrict to $A_i\in \cO_i$ then one
has $M_i\in \cC_i$. Also, as mentioned above, the map is equivariant under
diagonal conjugation and so there is a ``reduced Riemann--Hilbert map'':

\begin{equation} \label{eq: rRH}
\bcO:=\cO_1\times\cdots \times \cO_m\spq G
\quad\mapright{\nu_{\bf a}}\quad
\cC_1\fus\cdots\fus\cC_m\spq G =: \bcC
\end{equation}

\noindent
where the space $\bcO$ is the quotient of 
$
 \bigl\{ (A_1,\ldots,A_m) \ \bigl\vert \ 
          A_i\in\cO_i, \text{$\sum A_i= 0$} \bigr\}$
by overall conjugation by $G$ and
the space $\bcC$  is the quotient of
$\bigl\{ (M_1,\ldots,M_m) \ \bigl\vert \ 
          M_i\in \cC_i, M_m\cdots M_1 = 1 \bigr\}$
by overall conjugation by $G$. 
Generally this map $\nu_{\bf a}$ is an
injective holomorphic symplectic map between complex symplectic manifolds of
the same dimension.

The simplest case is when the representation is rigid, i.e. when the expected
dimensions of both sides of \eqref{eq: rRH} is zero. 
Then one knows the RHS of \eqref{eq: rRH} consists of precisely one point and
the LHS (at most) one point.

Our basic strategy is to look at the next simplest case, with the aim of
degenerating into the rigid case. Since the
spaces are symplectic, this corresponds to complex dimension two, 
i.e. both sides of \eqref{eq: rRH} are complex surfaces.

The principal example of  such ``minimally non-rigid'' systems occurs if we 
look at rank two systems with four poles on the sphere (i.e. $n=2,m=4$).
Without loss of generality (by tensoring by logarithmic connections on
line-bundles) 
one can work with $G=\SL_2(\IC)$ rather than
$\GL_2(\IC)$ and, using automorphisms of the sphere we can fix three of the
poles at $0,1,\infty$ and label the remaining pole position $t$.
Thus we are considering systems of  the form:
\beq\label{eq: lin syst}
\frac{d}{dz}-\left(
\frac{A_1}{z}+\frac{A_2}{z-t}+\frac{A_3}{z-1}\right),\qquad
A_i\in\g:=\lsl_2(\IC)
\eeq
By convention we denote the eigenvalues of $A_i$ by
$\pm\theta_i/2$ for $i=1,2,3,4$.
Schlesinger's equations imply that the residue $A_4=-\sum_1^3A_i$ at infinity
remains fixed; we will conjugate the system so that 
$A_4=\frac{1}{2}\,\diag(\theta_4,-\theta_4)$. The remaining conjugation
freedom is then just conjugation by the one-dimensional torus 
$T:=\diag(a,1/a), a\in\IC^*$; the space of such systems is then three
dimensional (quotienting by $T$ yields the surface $\bcO$).

Following \cite{JM81} pp.443-446 one may choose certain 
coordinates
$x,y,k$ 
on this space of systems and write down what Schlesinger's
equations become. 
One obtains a pair of coupled first-order nonlinear differential 
equations in $x,y$ (not dependent on $k$) 
and an equation for $k$ of the form $\frac{dk}{dt}=f(y,t)k$.
The coordinate $k$ corresponds to the torus action, which we can 
forget about since we are happy to consider Fuchsian systems
up to conjugation.
Eliminating $x$ from the coupled system yields the sixth Painlev\'e equation:
$$\frac{d^2y}{dt^2}=
\frac{1}{2}\left(\frac{1}{y}+\frac{1}{y-1}+\frac{1}{y-t}\right)
\left(\frac{dy}{dt}\right)^2
-\left(\frac{1}{t}+\frac{1}{t-1}+\frac{1}{y-t}\right)\frac{dy}{dt}  $$
$$
\quad\ +\frac{y(y-1)(y-t)}{t^2(t-1)^2}\left(
\al+\be\frac{t}{y^2} + \gamma\frac{(t-1)}{(y-1)^2}+
\delta\frac{ t(t-1)}{(y-t)^2}\right) $$
where the constants $\al,\be,\ga,\de$ 
are related to the $\theta$-parameters as follows:
\begin{equation}\label{thal}
\al=(\theta_4-1)^2/2, \qquad \be=-\theta_1^2/2, \qquad 
\ga=\theta^2_3/2, \qquad \delta=(1-\theta_2^2)/2.
\end{equation}

Since we will want to go back from a solution of \PVI\ to an explicit 
isomonodromic family of Fuchsian systems, we will give the explicit formulae
for the matrix entries of the system in terms of $y,y'$, in appendix
\ref{apx: residues}.

Now the bad news is that most solutions to \PVI\ cannot be written in terms of
classical special functions. From Watanabe's work \cite{watanabePVI} 
one knows that
either a solution is non-classical or it is a Riccati solution (corresponding
to a reducible or rigid 
monodromy representation $\rho$) or the solution $y(t)$ is an
algebraic function.

Since we are interested in explicit solutions corresponding to
irreducible non-rigid 
representations, the only possibility is to seek algebraic
solutions to \PVI, in other words solutions defined implicitly by 
equations of the form
$$F(y,t)=0$$
for polynomials $F$ in two variables.
We can rephrase this more geometrically:

\begin{defn}
An algebraic solution of \PVI\ consists of a triple $(\Pi,y,t)$ where
$\Pi$ is a compact (possibly singular) algebraic curve and  $y,t$ are 
rational functions on $\Pi$ such that:

$\bullet$
$t:\Pi\to \IP^1$ is a Belyi map (i.e. $t$ expresses $\Pi$ as a branched cover
of $\IP^1$ which only ramifies over $0,1,\infty$), and

$\bullet$
Using $t$ as a local coordinate on $\Pi$ away from ramification points,
$y(t)$  should solve \PVI, for some value of the parameters 
$\al,\be,\ga,\de$.
\end{defn}

Indeed given an algebraic solution in the form $F(y,t)=0$ one may take 
$\Pi$ to be the closure in $\IP^2$ of the affine plane curve defined by $F$.
That $t$ is a Belyi map on $\Pi$ follows from the Painlev\'e property of \PVI:
solutions will only branch at $t=0,1,\infty$ and all other singularities
are just poles. The reason we prefer this reformulation is that often the
polynomial $F$ is quite complicated and parameterisations 
of the plane curve defined by $F$ are usually simpler to write down. 
(The polynomial $F$ can be recovered as the minimal polynomial of $y$ over
$\IC(t)$, since $\IC(y,t)$ is a finite extension of $\IC(t)$.)

We will say the solution curve $\Pi$ 
is `minimal' or an `efficient parameterisation' if $y$ generates
the field of rational functions on $\Pi$, over $\IC(t)$, so that $y$ and $t$
are not pulled back from another curve covered by $\Pi$ (i.e. that
$\Pi$ is birational to the curve defined by $F$).

The main invariants of an algebraic solution are the genus of the (minimal) 
Painlev\'e
curve $\Pi$ and the
degree of the corresponding
Belyi map $t$ (the number of branches the solution has over the
$t$-line). 

Now the basic question is: what representations $\rho$ can we start with in
order to obtain an algebraic solution to \PVI?
Well, the solution must have only a finite number of branches and so we can
start by looking for finite branching solutions, 
and hope to prove in each case that the solution is
actually algebraic.

The important point is that one can read off the branching of the solution $y$
as $t$ moves around loops in the three-punctured sphere
$\IP^1\setminus\{0,1,\infty\}$ in terms of the corresponding linear
  representations $\rho$. 
One finds (cf. e.g. \cite{icosa} section 4) 
that $\rho$ transforms according to the natural 
action of the pure mapping
class group (which is isomorphic 
to $\pi_1(\IP^1\setminus\{0,1,\infty\})$ and  thus to 
the free group on two-letters $\cF_2$).  
Explicitly the generators $w_1,w_2$  of $\cF_2$ 
act on the monodromy matrices $\bM$ via $w_i=\omega_i^2$ where
$\om_i$ fixes $M_j$ for $j\ne i,i+1, \ (1\leqslant j\leqslant 4)$ and
\beq \label{eq: braiding}
\om_i(M_i,\ M_{i+1})= (M_{i+1},\ M_{i+1}M_{i}M^{-1}_{i+1}).
\eeq

\noindent
(Incidentally the geometric origins of this in the context of \PVI\ 
can be traced back at least to
Malgrange's work \cite{Mal-imd1long}
on the global properties of the Schlesinger equations.) 
The full classification of the representations $\rho$ living in finite orbits
of this action is still open, but there are some obvious ones: namely if 
$\rho$ takes values in a finite subgroup of $\SL_2(\IC)$ then the $\F_2$ 
orbit will clearly be finite.

Thus the program is to take such a finite subgroup $\Ga\subset G$, 
and go through the possible
representations $\rho:\pi_1(\IP^*)\to \Ga$ (whose image generates $\Ga$ say)
and find the corresponding 
\PVI\ solutions.
The two main problems to overcome in completing this program are:

1) There are lots of such representations (even up to conjugation), for example
for the binary 
tetrahedral group from \cite{Hall36} one knows 
there are  $520$  conjugacy classes 
of triples of generators.

2) We still need to find the \PVI\ solution explicitly.

For 1) we proceed as in \cite{icosa}; 
by using Okamoto's affine $F_4$ symmetry group of \PVI\  we can
drastically reduce the number of classes that arise for each group.
It is worth emphasising that upon applying an Okamoto transformation 
the monodromy
group may well become infinite, and currently there are very few examples
of algebraic
solutions to \PVI\  which are not equivalent to 
(or simple deformations of) a solution with
finite linear monodromy group (see the final remark of section \ref{sn: img} 
below).

For 2) we use Jimbo's asymptotic formula 
(see \cite{Jimbo82} and the corrected version in \cite{k2p} Theorem 4).
By looking at the degeneration
of the Fuchsian system into systems with only three poles 
(hypergeometric systems) and using explicit solutions of their Riemann--Hilbert
problems, Jimbo found explicit formulae for the leading term in the 
asymptotic expansion of \PVI\ solutions at zero.
Using the $\text{P}_{\text{VI}}$ equation 
these leading terms determine the Puiseux expansions of 
each branch of the
solution at zero and, taking sufficiently many terms, 
these enable us to find the solution completely if it is algebraic.

Philosophically the author views this work as an illustration of the utility
of Jimbo's asymptotic formula.
An alternative method of constructing solutions of Painlev\'e VI has
been proposed by Kitaev (and Andreev) 
\cite{Kit-sfit6, And-Kit-CMP, Kitaev-dessins} who call it the ``RS'' method
(see also Doran \cite{Chuck1} for similar ideas, section \ref{sn: img} below
and also \cite{Klein-ico} for
closely related ideas of F. Klein).
Kitaev \cite{Kit-sfit6} 
conjectures that all algebraic solutions arise in this way and, 
with Andreev, 
has found some solutions essentially by starting to enumerate all suitable 
rational maps along which a hypergeometric system may be pulled back.

One of our original aims was to 
try to ascertain what algebraic \PVI\ solutions are
known, up to equivalence under Okamoto transformations
and simple deformation (cf. e.g. \cite{icosa} Remark 15).
In other words the aim was
to see how much is known of what might be called the
 ``non-abelian Schwarz list'', viewing \PVI\ as the simplest
non-abelian Gauss--Manin connection.
The result is that, so far, all the algebraic solutions the author has seen 
have turned out to be related to a finite subgroup of $\SL_2(\IC)$
or to the $237$ triangle group (see section \ref{sn: img} below).\footnote{
One might be so bold as to conjecture that there are no others, simply because
no others have yet been seen, in spite of the variety of approaches used.}
As an illustrative example of what can happen consider 
solution 4.1.7.B of \cite{And-Kit-CMP}:
At first glance we see $t$ is a degree $8$ function of the parameter $s$ and
so one imagines a solution with $8$ branches (and wonders if it is related
to one of the eight-branch solutions of \cite{icosa} 
or \cite{Hit-Poncelet}
or of section \ref{sn: oct} below).
However one easily confirms that in fact
$$y=y_{21}=t+\frac{3\eta_\infty\sqrt{t(t-1)}}{\eta_\infty+1}$$
so it really only has two branches (it was inefficiently parameterised). 
In turn one finds (for any value of the constant $\eta_\infty$) this is 
equivalent to the well-known solution $y=\sqrt{t}$.

On the other hand Jimbo's formula gives us great control, in that we can often 
go directly from a
linear representation $\rho$ to the corresponding \PVI\ solution.
In particular the mapping class group orbit of $\rho$ tells us a priori the
number of branches (and lots more) that the solution will have.
At some point the author realised (see the introduction to \cite{icosa}) that 
there should be more solutions related to the symmetries of the Platonic 
solids than had already appeared; we have found it to be more efficient to
first ascertain directly what solutions arise in this way, than for example 
to enumerate rational maps.
(The author's understanding is that a theorem of Klein implies
that the solutions of sections \ref{sn: tet} and 
\ref{sn: oct} below
 and of \cite{icosa} will arise via rational pullbacks of a
hypergeometric system, but it is not clear if the enumeration started in 
\cite{And-Kit-CMP}
would ever have 
found all the corresponding rational maps independently.)
\end{section}

\begin{section}{The tetrahedral solutions} \label{sn: tet}

In this section we will classify the solutions to \PVI\ having linear
monodromy group equal to the binary tetrahedral group 
$\Ga\subset G=\SL_2(\IC)$.
The procedure is similar to that used in \cite{icosa} for the icosahedral
group.

First we examine (as in \cite{icosa} section 2)
the set $S$ of $G$-conjugacy classes of triples of generators
$(M_1,M_2,M_3)$ of $\Ga$ (i.e. two triples are identified if they are
related by conjugating by an element of $G$).
(Equivalently this is 
the set of conjugacy classes of representations $\rho$ of the
fundamental group of the 
four-punctured sphere into $\Ga$, once we choose a suitable 
set of generators.)
From Hall's formulae \cite{Hall36} one knows there are $12480$ 
triples of generators of $\Ga$ and dividing by $24$ (the size of the
image in $\PSL_2(\IC)$ of the normaliser of $\Ga$ in $G$) we find that
$S$ has cardinality  $520$.
Then we quotient $S$ further by the relation of geometric equivalence
(cf. \cite{icosa} section 4):
two representations are identified if they are related by the full mapping
class group, or by the set of even sign changes of the four monodromy matrices
$M_i$ (with $M_4=(M_3M_2M_1)^{-1}$). One finds there are precisely six such
geometric equivalence classes, and by Lemma 9 of \cite{icosa} this implies
there are at most six solutions to \PVI\ with tetrahedral monodromy which are
inequivalent under Okamoto's affine $F_4$ action. 

On the other hand we can 
look at the set of $\theta$-parameters corresponding to
the representations in $S$.
Since Okamoto transformations act by the standard $W_a(F_4)$ action on the
space of parameters, it is easy to find the set of inequivalent parameters that
arise from $S$, cf. \cite{icosa} section 3. 
(Since they are real we can map them all into the closure of a
chosen alcove.) 
We find there are exactly six sets of inequivalent parameters
that arise and so there are at least six inequivalent tetrahedral solutions. 
Combining with the previous paragraph we thus see there are precisely six
inequivalent tetrahedral solutions to \PVI.

Various data about the six classes and the corresponding \PVI\ solutions
are listed in tables $1$ and $2$.
Table $2$ lists a representative set of $\theta$-parameters for each
class together with numbers $\si_{ij}$ which uniquely determine a 
triple $M_1,M_2,M_3$ in $S$ 
(and thus the linear representation $\rho$) 
for that class with the given $\theta$ values, 
via the formula
$$\tr(M_iM_j)=2\cos(\pi \si_{ij}).$$

\begin{table}[h]
\begin{center}
\begin{tabular}{|c|c|c|c|c|c|c|c|c|c| }
\hline
  & \text{Degree}  &  \text{Genus} & \text{Walls} & \text{Type}
  & \text{Alcove Point} & $n$ &  \text{Nonlinear Group} & \text{Partitions} 
\\ \hline
1 &  1 &  0 &  2 & $a{b}^{2}$ & 35, 15, 15, 5 & 96 &  1 &  \   \\ \hline
2 &  1 &  0 &  3 & ${b}^{3}$ & 30, 10, 10, 10 & 32 &  1 &  \  \\ \hline
3 &  2 &  0 &  3 & ${b}^{4}-$ & 50, 10, 10, 10 & 48 &  $S_2$ &  $1, 2$ \\ \hline
4 &  3 &  0 &  3 & ${b}^{4}+$ & 40, 0, 0, 0 & 72 &  $S_3$ &  $3, 2$ \\ \hline
5 &  4 &  0 &  2 & $a{b}^{3}$ & 45, 5, 5, 5 & 128 &  $A_4$ &  $3$ \\ \hline
6 &  6 &  0 &  3 & ${a}^{2}{b}^{2}$ & 50, 10, 0, 0 & 144 & $A_4$ &  $2^2, 3^2$  \\ \hline

\hline
\end{tabular}

\vspace{0.2cm}
\caption{Properties of the tetrahedral solutions.}
\label{terahedral solution table}
\end{center}\end{table}

\begin{table}[h]
\begin{center}
\begin{tabular}{|c|c|c| }
\hline
  & $(\theta_1,\theta_2,\theta_3,\theta_4)$  
  & $(\si_{12},\si_{23},\si_{13})$ 
\\ \hline
1 & 1/2, 0, 1/3, 1/3 & 1/2, 1/3, 1/3 \\ \hline
2 & 1/3, 0, 1/3, 1/3 & 1/3, 1/3, 1/3 \\ \hline
3 & 1/3, 2/3, 2/3, 2/3 & 1/2, 1/3, 1/2 \\ \hline
4 & 2/3, 1/3, 1/3, 2/3 & 1/2, 1/3, 1/2 \\ \hline
5 & 1/3, 1/3, 1/3, 1/2 & 1/3, 2/3, 1/3 \\ \hline
6 & 1/2, 1/3, 1/3, 1/2 & 1/3, 1/2, 1/3 \\ \hline

\hline
\end{tabular}

\vspace{0.2cm}
\caption{Representative parameters for the tetrahedral solutions}
\label{terahedral param table}
\end{center}\end{table}

The first two columns of table $1$ list the degree and genus of the \PVI\
solution.
The column labelled ``Walls'' lists the 
number of affine $F_4$ reflection hyperplanes the parameters of the solution
lie on. The type of the solution enables us to see at a glance which class a
given element of $S$ lies in: Given $M_1,M_2,M_3,M_4\in \Ga$ their images in 
$\PSL_2(\IC)\cong \SO_3(\IC)$ are real rotations and we write an 
``$a$'' for each rotation by half of a turn, 
a ``$b$'' for each rotation by a third of a turn, and write
nothing for each trivial rotation thus obtained. 
This distinguishes all classes except $3$ and $4$ which both 
correspond to four
rotations by a third of a turn: each $M_i$ thus has parameter
$\theta_i=1/3$ or $\theta_i=2/3$. For class $3$ there are always an odd
number of each type of $\theta$ 
($1/3$ or $2/3$) so we write a minus, and for class $4$
there are always an even number of each type, so we write a plus.

Finally the rest of table $1$ lists the corresponding alcove point (scaled by
$60$), 
the number $n$ of elements of $S$ belonging to each class,
the monodromy group of the cover $t:\Pi\to \IP^1$ 
and the unordered collection of 
sets of ramification indices of this cover over $t=0,1,\infty$ 
(repeating the last set of indices until three are obtained). 
Thus for example each solution corresponding to 
row $6$ has indices $(3,3)$ over two points amongst $\{0,1,\infty\}$
and indices $(1,1,2,2)$ over the third.

All of the tetrahedral solutions have genus zero so we may
 take $\Pi$ to be $\IP^1$ with
parameter $s$ and write the solutions as functions of $s$.
As in the icosahedral case the solutions with at most $4$ branches are closely
related to known solutions.
For classes $1$ and 
$2$ one of the monodromy matrices is projectively trivial and so
these rows correspond to pairs of generators of the tetrahedral group, i.e. to
the two tetrahedral entries on Schwarz's list of algebraic hypergeometric
functions. The corresponding \PVI\ solutions are both just 
$y=t$ with the parameters as listed in table $2$.
As in \cite{icosa} one finds class $3$ contains the solution 
$y=\pm\sqrt{t}$ (with the parameters as listed in table $2$).
Class $4$ contains the tetrahedral solution 
\begin{eqn}\label{eq: hittet}
y=\frac{(s-1)(s+2)}{s(s+1)},\qquad
t=\frac{(s-1)^2(s+2)}{(s+1)^2(s-2)}
\end{eqn}
\!\!on p.592 of \cite{Hit-Octa}
(with the parameters as listed in table $2$)
and is equivalent to a solution found independently 
by Dubrovin \cite{Dub95long} (E.31).
Also class $5$ contains a simple deformation of the 
four-branch dihedral solution in section 6.1 of \cite{Hit-Poncelet}:
\begin{equation} \label{eq: dih soln}
y=\frac{s^2(s+2)}{s^2+s+1},\qquad
t=\frac{s^3(s+2)}{2s+1},
\end{equation}
that is, this solution is tetrahedral if we use the parameters in table $2$,
rather than the parameters $1/2,1/2,1/2,1/2$ for which it is dihedral.

Thus we are left with one solution, corresponding to row $6$.
Using Jimbo's asymptotic formula to compute the Puiseux expansions etc. 
(as in \cite{k2p} section $5$, especially p.193) 
we find the following  solution in this class:

\begin{gather}
\text{Tetrahedral solution $6$,  $6$ branches
$(\theta_1,\theta_2,\theta_3,\theta_4)=(1/2, 1/3, 1/3, 1/2)$:}\notag \\
y=
-{\frac {s \left( s+1 \right)  \left( s-3 \right) ^{2}}
{3 \left( s+3 \right)  \left( s-1 \right) ^{2}}},\qquad
t=-\left(
{\frac { \left( s+1 \right)  \left( s-3 \right)}
{ \left( s-1 \right) \left( s+3 \right)        }}\right)^3\notag
\end{gather}
(We have recently learnt that this is equivalent to solution 4.1.1A 
in \cite{And-Kit-CMP}.) 
It is now easy to write down the
explicit isomonodromic family of Fuchsian systems in this
case, thereby solving the Riemann--Hilbert problem for this class of
representations $\rho$, for an arbitrary configuration of the four 
pole positions (up to automorphisms of $\IP^1)$. 
(We will leave for the reader the
analogous substitutions for the other solutions below.)
Using the formulae in appendix \ref{apx: residues}
one finds the family of systems
parameterised by $s\in\IP^1$ is (up to overall conjugation): 
$$\frac{d}{dz}-\left(\frac{A_1}{z}+\frac{A_2}{z-t(s)}+\frac{A_3}{z-1}\right)$$
where
$$A_1=
\left( \begin {smallmatrix} 
\left( {s}^{2}+3 \right)  \left( {s}^{6}-51{s}^{4}+99{s}^{2}-81 \right) 
&
4s \left( {s}^{4}-9 \right) \\ 
4\left( 5{s}^{6}-75{s}^{4}+135{s}^{2}-81 \right) s 
\left( {s}^{4}-9 \right) 
& 
\ -\left( {s}^{2}+3 \right)  \left( {s}^{6}-51{s}^{4}+99{s}^{2}-81 \right) 
\end {smallmatrix} \right)\Big/\Delta$$
$$A_2=
\left( \begin {smallmatrix} 
4 \left( s+3 \right)  \left( s-1 \right) ^{2}{s}^{2} \left(
{s}^{3}-{s}^{2}+3s+9 \right) 
&
\ -2 \left( s+3 \right)  \left( s-1 \right) ^{2} \left( {s}^{2}+2s+3 \right) \\
 -2\left( s+3 \right)  \left( s-1 \right) ^{2} \left( {s}^{3}-3{s}^{2}-9s-9 \right)  \left( 5{s}^{5}-5{s}^{4}-45
s-27 \right) 
&
\ -4 s^2\left( s+3 \right)  \left( s-1 \right) ^{2}
\left( {s}^{3}-{s}^{2}+3s+9 \right)
\end {smallmatrix} \right)\Big/\Delta$$
$$A_3=\diag(-\theta_4,\theta_4)/2-A_1-A_2,\qquad
\Delta=-36 \left( {s}^{2}+3 \right)  
\left( s^2-1 \right) ^{2} \left( {s}^{2}-9 \right).$$

Note that if the denominator $\Delta$ is zero then  $t\in \{0,1,\infty\}$
since 
$$1-t=2\,{\frac { \left( {s}^{2}+3 \right) ^{2} \left( {s}^{2}-3 \right) }{ \left( s+3 \right) ^{3} \left( s-1 \right) ^{3}}}.$$
Thus the system 
is well defined for all $s$ in
$t^{-1}(\IP^1\setminus\{0,1,\infty\})$ except possibly at $s=\infty$ (where
$t=-1$). However writing $s=1/s'$ it is easy  to conjugate the system to one
well-defined also at $s=\infty$. Thus one never encounters configurations
requiring a nontrivial bundle; the Malgrange divisor is trivial in this
situation (in spite of the fact the solution 
$y$ does have a pole at $s=\infty$); indeed one knows
the corresponding $\tau$ function (whose zeros lying over 
$\IP^1\setminus\{0,1,\infty\}$ correspond to nontrivial bundles)
satisfies:
$$d\log(\tau)=
\tr\left(A_2\left(\frac{A_1}{t}+\frac{A_3}{t-1}\right)\right)dt=
-{\frac {{s}^{6}+6\,{s}^{5}+3\,{s}^{4}-8\,{s}^{3}-9\,{s}^{2}-54\,
s-27}
{3\, \left( {s}^{4}-9 \right)  \left( {s}^{2}-1 \right)  \left( {s}^
{2}-9 \right) }}ds $$
which is nonsingular for all $s\in t^{-1}(\IP^1\setminus\{0,1,\infty\})$.

\begin{rmk}
Sometimes one is interested in Fuchsian equations with given monodromy,
rather than systems. 
To obtain these one may choose a cyclic vector, or more simply substitute the
\PVI\ solution into the standard formulae for the isomonodromic family of
Fuchsian equations. In the present case one obtains the equation:
$$\frac{d^2}{dz^2} + a_1\frac{d}{dz} + a_2$$
where $a_1,a_2$ are respectively:
$$
\frac{1}{2z}+\frac{2}{3\left( z-1 \right)}
+\frac{2}{ 3(z-t(s))}
-\frac{1}{z-y(s)},
$$
$$
\frac{1}{z(z-1)}
\left(
\frac{7s^6-6s^5+3s^4+4s^3-63s^2-54s-27}
{18(s+3)^3(s-1)^3(z-t(s))}
-\frac{(s^2-2s+3)s^2}{9(s+3)(s-1)^2(z-y(s))}-\frac{1}{18}
\right).
$$
For generic 
$s$ this is a Fuchsian equation with non-apparent singularities at 
$z=0,1,t,\infty$ and an apparent singularity at $z=y$, realising the given 
(projective) monodromy representation.
In special cases (when $y=0,1,t,\infty$) it will have just the four
non-apparent singularities (and will thus be a so-called ``Heun equation'').
For example specialising to
$s=0$ one finds $y=0,t=-1$ and the equation becomes that with
$$a_1=
-\frac{1}{2z}+\frac{2}{3\left( z-1 \right)}
+\frac{2}{ 3(z+1)}\qquad
a_2=-\,{\frac {1}{ 18\left( z-1 \right)  \left( z+1 \right) }}$$
which is a Heun equation whose projective monodromy representation
is that specified by 
row $6$ of table $2$.

\end{rmk}

\begin{rmk}
At the editor's request we will explain how one may verify directly that these
\PVI\ solutions actually do correspond to Fuchsian systems with linear
monodromy representations as specified by table 2.
For the rigid cases, rows $1$ and $2$, this is immediate, by rigidity.
For the others, first one may check that the solutions actually do solve \PVI.
This can be done directly (by computing 
the derivatives of $y$ with respect to $t$
and substituting into the \PVI\ equation).\footnote{    \label{fn: accfile}
To aid the reader interested in examining the solutions of this article
(and to help avoid typographical errors)
a Maple text file of the solutions has been included with the source file of
the preprint version on
the math arxiv (math.DG/0501464). 
This may be downloaded by clicking on ``Other formats'' and
unpacked with the commands `gunzip 0501464.tar' and `tar -xvf 0501464.tar', 
at least on a Unix system.
}
Having done this we know the formulae of appendix \ref{apx: residues} 
do indeed give
an isomonodromic family of Fuchsian systems.
To see it has the monodromy representation specified by table $2$
we first compute the Puiseux expansions at $0$ of each branch of the function 
$y(t)$ (only the leading terms will be needed).
On the other hand Jimbo's asymptotic formula (in the form in 
\cite{k2p} Theorem 4)
computes the leading term in the asymptotic expansion of the 
\PVI\  solution corresponding to the given monodromy representation $\rho$
(the leading term is of the form $a t^b$ where $a$ and $b$
are explicit functions of 
$\theta_1,\theta_2,\theta_3,\theta_4,\sigma_{12},\sigma_{23},\sigma_{13}$).
Then it is sufficient to check this leading term equals one of the
leading terms of the Puiseux expansions of $y(t)$.
The logic is that, in the cases at hand, the leading term determines the whole
Puiseux expansion (using the recursion {determined} by the  \PVI\ equation)
and this is convergent so determines the solution locally, and thus globally
by analytic continuation. (For the solutions we construct here this is
automatic since we constructed the solution starting with the results of
Jimbo's formula.)

Some simpler, but not entirely conclusive, checks are as follows:

1) compare the monodromy of the Belyi map $t$ with the $\F_2$ action 
\eqref{eq: braiding} on the conjugacy class of the representation $\rho$
(if we didn't know better it would appear as a
miracle that the solution, constructed out
of just the Puiseux expansion at $0$, turns out to have the right branching at
$1$ and $\infty$ too). 

2) compute directly the Galois group 
of one of the Fuchsian systems in the
isomonodromic family.
(Together with the exponents this goes a long way to pinning down the
monodromy representation.)
There are various ways to do this, one of which is to convert the system into
an equation (e.g. via a cyclic vector) and use the facility on 
Manuel Bronstein's webpage:

$\text{http:/\!/www-sop.inria.fr/cafe/Manuel.Bronstein/sumit/bernina\_demo.html}$

\noindent
(This requires finding a suitable 
rational point on the Painlev\'e curve, which, if
possible, is easy
in the genus zero cases, and not too difficult using Magma in the genus one
cases.)

\end{rmk}

\end{section}

\begin{section}{The octahedral solutions} \label{sn: oct}

For the octahedral group we do better and find more new solutions.
In this case, by \cite{Hall36} or direct computation, 
$S$ has size $3360$, which reduces to just thirteen classes 
under either geometric or parameter equivalence. Thus there are exactly 
thirteen
octahedral solutions to \PVI, up to equivalence under Okamoto's affine $F_4$
action.

Data about these classes are listed in tables $3$ and $4$.
In this case the type of the solution may contain the symbol ``$g$'' which
indicates that one of the corresponding rotations in $\SO_3$ is a rotation by a
quarter of a turn. Also, in some cases rather than list the monodromy
group of the cover $t:\Pi\to \IP^1$ we just give its size.

\begin{table}[h]
\begin{center}
\begin{tabular}{|c|c|c|c|c|c|c|c|c|c| }
\hline
  & \text{Degree}  &  \text{Genus} & \text{Walls} & \text{Type}
  & \text{Alcove Point} & $n$ &  \text{Group (size)} & \text{Partitions} 
\\ \hline
1 & 1 & 0 & 1 & $abg$ & (65, 35, 25, 5)/2 & 192 & $1$ & \  \\ \hline
2 & 1 & 0 & 2 & $bg^2$ & 25, 10, 10, 5 & 96 & $1$ & \  \\ \hline
3 & 2 & 0 & 2 & $b^2g^2$ & 45, 15, 10, 10 & 96 & $S_2$ & $1, 2$ \\ \hline
4 & 3 & 0 & 1 & $abg^2$ & 40, 10, 5, 5 & 288 & $S_3$ & $3, 2$\\ \hline
5 & 4 & 0 & 2 & $ag^3$ & (75, 15, 15, 15)/2 & 128 & $A_4$ &  $3$ \\ \hline
6 & 4 & 0 & 3 & $g^4$ & 30, 0, 0, 0 & 32 & $A_4$ & $3$ \\ \hline
7 & 6 & 0 & 2 & $a^2bg$ & (95, 25, 5, 5)/2 & 576 & 24 & 
$2^2, 3^2, 2\, 4$ \\ \hline
8 & 6 & 0 & 2 & $b^2g^2$ & 35, 5, 0, 0 & 288 & 36 & 
$3, 2\, 4$ \\ \hline
9 & 8 & 0 & 1 & $ab^2g$ & (85, 15, 15, 5)/2 & 768 & 576 &
$2^2\,3, 2^2\,4$ \\ \hline
10 & 8 & 0 & 3 & $a^2g^2$ & 45, 15, 0, 0 & 192 & 192 & 
$3^2, 2\,3^2$ \\ \hline
11 & 12 & 0 & 3 & $a^2b^2$ & 50, 10, 0, 0 & 288 & 576 & 
$2^2\,3^2, 2^2\,4^2$ \\ \hline
12 & 12 & 1 & 3 & $a^3b$ & 55, 5, 5, 5 & 288 & 96 & 
$3^4, 2^2\,4^2$ \\ \hline
13 & 16 & 0 & 3 & $a^3g$ & (105, 15, 15, 15)/2 & 128 & 3072 & 
$2^2\, 3^4$ \\ \hline

\hline
\end{tabular}

\vspace{0.2cm}
\caption{Properties of the octahedral solutions.}
\label{octahedral solution table}
\end{center}\end{table}

\begin{table}[h]
\begin{center}
\begin{tabular}{|c|c|c| }
\hline
  & $(\theta_1,\theta_2,\theta_3,\theta_4)$  
  & $(\si_{12},\si_{23},\si_{13})$ 
\\ \hline
1 & 1/2, 0, 1/3, 1/4 & 1/2, 1/3, 1/4 \\ \hline
2 & 1/3, 0, 1/4, 1/4 & 1/3, 1/4, 1/4 \\ \hline
3 & 1/3, 1/4, 1/4, 2/3 & 1/2, 1/2, 1/2 \\ \hline
4 & 1/2, 1/4, 1/4, 2/3 & 1/2, 1/3, 3/4 \\ \hline
5 & 1/4, 1/4, 1/4, 1/2 & 1/3, 1/2, 1/3 \\ \hline
6 & 1/4, 1/4, 1/4, 1/4 & 1/3, 0, 1/3 \\ \hline  
7 & 1/2, 1/2, 1/4, 2/3 & 1/2, 1/2, 1/3 \\ \hline
8 & 1/3, 3/4, 1/3, 3/4 & 1/2, 3/4, 1/3 \\ \hline
9 & 1/3, 1/4, 1/2, 2/3 & 1/2, 2/3, 3/4 \\ \hline
10 & 1/2, 1/4, 1/2, 3/4 & 2/3, 2/3, 1 \\ \hline
11 & 1/3, 1/2, 1/2, 2/3 & 1/2, 1/2, 1/4 \\ \hline
12 & 1/2, 1/2, 1/2, 2/3 & 1/2, 1/4, 2/3 \\ \hline
13 & 1/2, 1/2, 1/2, 3/4 & 1/2, 2/3, 1/3 \\ \hline

\hline
\end{tabular}
\vspace{0.2cm}
\caption{Representative parameters for the octahedral solutions}
\label{octahedral param table}
\end{center}\end{table}

The octahedral 
solutions with at most $4$ branches correspond to the following 
known solutions.
As in \cite{icosa} one finds:
The first two classes correspond to the octahedral entries on Schwarz's list
of algebraic 
hypergeometric functions (and the \PVI\ solution is $y=t$ with the
parameters indicated in table $4$).
Solution $3$ is $y=\pm\sqrt{t}$ with the parameters listed 
in table $4$,
solution $4$ has $3$ branches and is a simple deformation of the $3$-branch
tetrahedral solution above (namely it is the solution in equation 
\eqref{eq: hittet}, but with the parameters given in table $4$),
solution $5$ is a simple deformation of the $4$-branch
dihedral solution (namely it is the solution in equation 
\eqref{eq: dih soln}, but with the parameters given in table $4$),
and solution $6$ is the $4$-branch octahedral solution 
\begin{eqn}\label{eq: hitoct}
y=\frac{(s-1)^2}{s(s-2)},\qquad
t=\frac{(s+1)(s-1)^3}{s^3(s-2)}
\end{eqn}
\!\!on p.588 of \cite{Hit-Octa}, with the parameters as in table $4$, which is
equivalent to a solution found independently by Dubrovin \cite{Dub95long}
(E.29).

For the remaining $7$ solutions, rows $7$--$13$, we will construct an explicit
solution in each class using
Jimbo's asymptotic formula.
More computational details appear in appendix \ref{apx: cdets}.
(We have recently learnt that solutions 8 and 10 are equivalent to those of 
\cite{Kitaev-dessins}, 3.3.3 top of p.22, and 3.3.5 bottom of p.23,
respectively.)
The formulae obtained are as follows:

\begin{gather}
\text{Octahedral solution $7$,  $6$ branches
$(\theta_1,\theta_2,\theta_3,\theta_4)=(1/2, 1/2, 1/4, 2/3)$:}\notag \\
y={\frac { 9 s \left( 2\,{s}^{3}-3\,s+4 \right)}
{ 4 \left( s+1 \right)  \left( s-1 \right) ^{2} \left( 2\,{s}^{2}+6
\,s+1 \right)  }},\qquad
t={\frac {27{s}^{2}}{ 4 \left( {s}^{2}-1 \right) ^{3}}}\notag
\end{gather}

\begin{gather}
\text{Octahedral solution $8$,  $6$ branches
$(\theta_1,\theta_2,\theta_3,\theta_4)=(1/3, 3/4, 1/3, 3/4)$:}\notag \\
y={\frac { \left( 2\,{s}^{2}-1 \right)  \left( 3\,s-1 \right) }{
 2 s \left( 2\,{s}^{2}+2\,s-1 \right)  \left( s-1 \right)}},\qquad
t=-{\frac { \left( 3\,s-1 \right) ^{2}}{ 8 \left( 2\,{s}^{2}+2\,s-1
 \right)  \left( s-1 \right) {s}^{3}}}\notag
\end{gather}

\begin{gather}
\text{Octahedral solution $9$,  $8$ branches
$(\theta_1,\theta_2,\theta_3,\theta_4)=(1/3, 1/4, 1/2, 2/3)$:}\notag \\
y=
{\frac {{s}^{3} \left( 2\,{s}^{2}-4\,s+3 \right)  \left( {s}^{2}-2\,s+
2 \right) }{ \left( 2\,{s}^{2}-2\,s+1 \right)  \left( 3\,{s}^{2}-4\,s+
2 \right) }},\qquad
t=\left(
{\frac {{s}^{2} \left( 2\,{s}^{2}-4\,s+3 \right)}{ 3\,{s}^
{2}-4\,s+2 }}\right)^2 \notag
\end{gather}

\begin{gather}
\text{Octahedral solution $10$,  $8$ branches
$(\theta_1,\theta_2,\theta_3,\theta_4)=(1/2, 1/4, 1/2, 3/4)$:}\notag \\
y={\frac {32 s \left( s+1 \right)  \left( 5\,{s}^{2}+6\,s-3 \right) }{
 \left( {s}^{2}+2\,s+5 \right)  \left( 3\,{s}^{2}+2\,s+3 \right) ^{2}}
},\qquad
t={\frac {1024 {s}^{3} \left( s+1 \right) ^{2}}{ \left( {s}^{2}+6\,s+1
 \right)  \left( 3\,{s}^{2}+2\,s+3 \right) ^{3}}}\notag
\end{gather}

\begin{gather}
\text{Octahedral solution $11$,  $12$ branches
$(\theta_1,\theta_2,\theta_3,\theta_4)=(1/3, 1/2, 1/2, 2/3)$:}\notag \\
y=\frac { \left( s+1 \right) 
 \left( 7\,{s}^{4}+16\,{s}^{3}+4\,{s}^{2}-4 \right)\, r }
{s^3 \left( s-2 \right) 
\left( {s}^{4}-4\,{s}^{2}+32\,s-28 \right)  
},\qquad
t=\left(\frac { \left( s+1 \right) ^{2}\, r}
{  \left( s-2 \right) ^{2}{s}^{2}}\right)^2\notag
\end{gather}
where
$
r=
4 \left( 3\,{s}^{2}-4\,s+2 \right)/ \left( {s}^{2}+4\,s+6 \right)
$.

\

The next solution, number $12$, 
has genus one. In this case we take $\Pi$ to be the elliptic curve
$$u^2= \left( 2\,s+1 \right)  \left( 9\,{s}^{2}+2\,s+1 \right).$$
As functions on this curve the solution is:

\begin{gather}
\text{Octahedral solution $12$, genus one, $12$ branches
$(\theta_1,\theta_2,\theta_3,\theta_4)=(1/2, 1/2, 1/2, 2/3)$:}\notag \\
y= \frac{1}{2}+
{\frac 
{45\,{s}^{6}+20\,{s}^{5}+95\,{s}^{4}+92\,{s}^{3}+39\,{s}^{2}-3}
{ 4\,\left( 5\,{s}^{2}+1 \right)  \left( s+1 \right) ^{2}u}}\notag\\
t= \frac{1}{2}+
{\frac { s\left( 2\,s+1 \right) ^{2} \left( 27\,{s}^{4}+28\,{s}^{3}+26
\,{s}^{2}+12\,s+3 \right)}{\left( s+1 \right) ^{3}u^3}}\notag
\end{gather}

Finally the last solution, number $13$, has $16$ branches and genus zero. 
This is possibly the highest degree 
genus zero solution amongst all
algebraic solutions of \PVI.
It is also special since it has no real
branches. 
Thus necessarily the parameterisation of the solution is not defined over
$\IQ$ although the solution curve $F(y,t)=0$ itself has $\IQ$ coefficients,
as does the Belyi map $t$.

\begin{gather}
\text{Octahedral solution $13$,  $16$ branches,
$(\theta_1,\theta_2,\theta_3,\theta_4)=(1/2,1/2,1/2,3/4)$:}\notag\\
y=-{\frac {\left( 1+ i \right)  \left( s^2-1 \right) 
\left( {s}^{2}+2\,is+1 \right)  
 \left( {s}^{2}-2\,is+1 \right) ^{2} c}
{4\,s \left( {s}^{2}+i \right)\left( {s}^{2}-i \right) ^{2} 
\left( {s}^{2}+(1+i)s-i \right) d
 }},\quad 
t={\frac { \left( s^2-1 \right) ^{2} 
\left( {s}^{4}+6\,{s}^{2}+1 \right) ^{3}}
{32\, {s}^{2} \left( {s}^{4}+1 \right) ^{3}}}\notag
\end{gather}
where 
$$c=
{s}^{8}-(2-2i)\,{s}^{7}-(6+2i)\,{s}^{6}
+(10+2i)\,{s}^{5}
+4\,i{s}^{4}
+(10-2i)\,{s}^{3}
+(6-2i)\,{s}^{2}
-(2+2i)\,s-1 
$$
$$d=
 {s}^{6}
-(3+3i)\,{s}^{5}
+3\,i{s}^{4}
+(4-4i)\,{s}^{3}
+3\,{s}^{2}+
(3+3i)\,s+i
$$

\begin{rmk}
The author has  recently understood that an alternative way to 
construct some (but not all) of these tetrahedral and octahedral solutions 
would have been to use
the quadratic transformations of \cite{Kitaev-quad-p6}.
For example tetrahedral solution 6 could have been obtained from 
tetrahedral solution 4 in
this way (a fact that was apparently not noticed in \cite{And-Kit-CMP}).
It is debatable whether 
this would have been simpler for us than the direct method used here,
given what had already been done in \cite{k2p,icosa}.
(The quadratic transformations were crucial however to construct the higher
genus icosahedral solutions, cf. \cite{ipc}.) 
\end{rmk}

\end{section}

\begin{section}{Infinite monodromy groups} \label{sn: img}
In this final section 
we will give some examples of solutions corresponding to
non-rigid 
representations $\rho$ into some infinite subgroups of $G=\SL_2(\IC)$.
The point is that the method we are using to construct \PVI\ solutions 
should work provided that the $\F_2$ orbit of $\rho$ is finite, and above we
just 
used the finiteness of the image of $\rho$ as a convenient way of ensuring
this. 

Thus we are looking for representations $\rho$ having finite $\F_2$ orbits, or
more concretely, matrices $M_1,M_2,M_3\in G$ having 
finite orbit under the
action \eqref{eq: braiding}.
(Such an $\F_2$ orbit, on conjugacy classes of representations,
gives the permutation representation of the Belyi cover  
$t:\Pi\to \IP^1$ for the corresponding \PVI\ solution.
We would like to find interesting
$\F_2$ orbits in order to find interesting \PVI\ solutions.)
So far there appear to be four ways (apart from guessing) of finding
representations  $\rho$ having finite $\F_2$ orbits.

Firstly 
one can just set the parameters to be sufficiently irrational in one of 
the families of solutions. (For example $y=\sqrt{t}$ is a solution provided
$\theta_1+\theta_4=1,\theta_2=\theta_3$ amongst other possibilities.)

Secondly one can sometimes
apply an Okamoto transformation to a known solution and change
$\rho$ into a representation having infinite image. 
For example if we take the $16$ branch octahedral solution above and apply the
Okamoto 
transformation corresponding to the central node of the extended $D_4$ Dynkin
diagram, then we obtain a \PVI\ solution whose
corresponding linear representation has image equal to the 
$(2,3,8)$ triangle group.
To see this we recall \cite{IIS, k2p}
that Okamoto's affine $D_4$ action does not change the
quadratic functions $\tr(M_iM_j)=2\cos(\pi \si_{ij})$ 
of the monodromy data, only the $\theta$-parameters.
In this case the $16$ branch octahedral solution has data
$$\theta=(1/2, 1/2, 1/2, 3/4), \qquad\sigma=(1/2, 2/3, 1/3)$$
on one branch and the solution after applying the transformation has data
$$\theta=(3/8, 3/8, 3/8, 5/8), \qquad\sigma=(1/2, 2/3, 1/3).$$
One may show that the image in $\PSL_2(\IC)$ of the corresponding triple
$(M_1,M_2,M_3)$ generate a $(2,3,8)$ triangle group (see appendix 
\ref{apx: 237}).
The corresponding solution to \PVI\ is given by the formula
$$y_{238}(s)=y(s) + \frac{2-\sum_1^4\theta_i}{2\, p(y,y',t)}$$
where $y,t, \theta_i$ are as for the $16$ branch octahedral solution
and $p$ as in appendix \ref{apx: residues}. 
Explicitly one finds the solution is 
\begin{gather}
\text{$2,3,8$ solution, genus zero, $16$ branches,
$(\theta_1,\theta_2,\theta_3,\theta_4)=(3/8, 3/8, 3/8, 5/8)$:}\notag\\
y=
-\frac { \left( 1+i \right)\left( s^2-1 \right)    
\left( {s}^{2}+2\,is+1 \right) 
\left( {s}^{2}-2\,is+1 \right) ^{2}d'}
{ 8 s \left( {s}^{2}+i \right)\left( {s}^{2}-i\right) ^{2} d
}\notag
\end{gather}
with $t$ and $d(s)$ as for the $16$ branch octahedral solution, and
$d'(s)=\overline{d(\overline{s})}$.
In turn, via the formulae in appendix \ref{apx: residues}, 
this yields the explicit family of Fuchsian systems having projective
monodromy group the $(2,3,8)$ triangle group.


Thirdly, the idea of \cite{k2p} was to use a different realisation of \PVI\ 
as controlling isomonodromic deformations of certain $3\times 3$ systems. The
corresponding monodromy representations were subgroups of $\GL_3(\IC)$
generated by complex reflections, and again one will obtain  
finite branching solutions
by taking representations into a finite group generated by complex reflections.
Applying this to the Klein complex reflection group led to an algebraic
solution to \PVI\ with $7$ branches. 
Moreover \cite{k2p} described explicitly
how to go between this $3\times 3$ picture and the standard $2\times 2$
picture used here, both on the level of systems and monodromy data. The upshot
is that if we substitute the Klein solution into the formulae of appendix
\ref{apx: residues}
below, then we obtain an isomonodromic family of $2\times 2$
Fuchsian systems with
monodromy data on one branch given by:
\beq\label{eq: klein data}
\theta=(2/7, 2/7, 2/7, 4/7), \qquad\sigma=(1/2, 1/3, 1/2).
\eeq
This determines a representation into (a lift to $G$ of) the 
$(2,3,7)$ triangle
group (and one may show as in appendix \ref{apx: 237} 
its image is not a proper subgroup). 
Moreover it was proved in \cite{k2p} that this cannot be obtained by
Okamoto transformations from a representation into a finite 
subgroup of $\SL_2(\IC)$.

Fourthly one may obtain such representations by 
pulling back certain hypergeometric systems along certain rational maps
(cf. Doran \cite{Chuck1} and Kitaev \cite{Kitaev-dessins}). 
This is closely related 
Klein's theorem that all second order Fuchsian equations with finite monodromy
group are pull-backs of hypergeometric equations with finite monodromy. 
The basic idea is as follows.

Label two copies of $\IP^1$ by $u$ and $d$ (for upstairs and downstairs).
Choose four integers $n_0,n_1,n_\infty,N\ge 2$ and 
suppose we have an algebraic family of branched covers
$$\pi:\IP^1_u\to \IP^1_d$$ 
of degree $N$, parameterised by a curve $\Pi$ say, such that:

$1)$ $\pi$ only branches at four points $0,1,\infty$ and at a variable point 
$x\in \IP^1_d$,

$2)$ All but four of the ramification indices over $0,1,\infty$ divide 
$n_0,n_1,n_\infty$ respectively. 
In other words 
if $\{e_{i,j}\}$ are the ramification indices over $i=0,1,\infty$ then
as $j$ varies precisely four of the numbers
$$\frac{e_{0j}}{n_0},\qquad \frac{e_{1j}}{n_1},
\qquad \frac{e_{\infty j}}{n_\infty}$$
are not integers.
Let $t$ be the cross-ratio of the corresponding four ramification 
points of $\IP^1_u$, in some order, and so we have
a coordinate on $\IP^1_u$ such that these four points 
occur at $0,t,1,\infty$.

$3)$ $\pi$ has minimal ramification over $x$, i.e. $\pi$ ramifies at just one
point over $x$, with ramification index $2$.



The idea of \cite{Chuck1}, \cite{Kitaev-dessins} is to
take a hypergeometric system on $\IP^1_d$ with projective 
monodromy around $i$ equal to an $n_i$'th root of the identity, 
for $i=0,1,\infty$ and 
pull it back  along $\pi$.
One then 
obtains an isomonodromic family on $\IP^1_u$ with non-apparent
singularities at $0,t,1,\infty$ and possibly some 
apparent singularities at the other ramification points. 
All of the apparent singularities can be
removed, for example by applying suitable Schlesinger transformations,
yielding an isomonodromic family of systems of the desired form. 

In particular the problem of constructing algebraic solutions of \PVI\ now
largely 
becomes a purely algebraic problem 
about families of covers, although it is only
conjectural that all algebraic solutions arise in this way.

However the algebraic construction of such covers seems difficult. 
First one has
the topological problem of finding such covers, then one needs to find the 
full family of covers explicitly 
(this amounts to finding a parameterised solution of a large system of
algebraic equations, typically with one less equation than the number of
variables, so the solution is a curve). See \cite{Kitaev-dessins} for some
interesting examples however (but one should be aware that some of these
  solutions are equivalent to each other and to
  known solutions via Okamoto transformations).

Our perspective here is that just the topology of the cover is enough to
determine the monodromy of the Fuchsian equation on $\IP^1_u$ and we can then
apply our previous method \cite{k2p} 
to construct the explicit \PVI\ solution. 
In other words just one topological cover $\pi$ 
gives us the desired representation $\rho$ living in a finite $\F_2$ orbit.

To find some interesting topological covers we consider the list appearing in
Corollary 4.6 of \cite{Chuck1}. 
Here Doran classified the possible ramification 
indices of the cover $\pi$ in the cases where the monodromy group of the 
hypergeometric system downstairs is an arithmetic triangle group in 
$\SL_2(\IR)$.
(Contrary to the wording in \cite{Chuck1} this does not determine the topology
of the cover.)
We will content ourselves with looking at the last four entries of Doran's
list, which say that the integers $N, n_0,n_1,n_\infty$ and the ramification
indices are:
\begin{align*}
&10,\  2,3,7 \qquad [2,\ldots,2], [3,3,3,1], [7,1,1,1]\\
&12,\  2,3,7 \qquad [2,\ldots,2], [3,3,3,3], [7,2,1,1,1]\\
&12,\  2,3,8 \qquad [2,\ldots,2], [3,3,3,3], [8,1,1,1,1]\\
&18,\  2,3,7 \qquad [2,\ldots,2], [3,\ldots,3], [7,7,1,1,1,1].\\
\end{align*}

The basic problem now is to find such covers topologically, in other words to
find the possible permutation representations.
(The cover of the four punctured sphere 
$\IP^1_d\setminus\{x,0,1,\infty\}$ is determined by its monodromy, which
amounts to
four elements of $\Sym_N$ having product equal to the identity and whose
conjugacy classes---i.e. cycle types---are 
as specified by the given ramification indices.)

The simplest way to do this is to draw a picture. Suppose we fix $x=-1$ and
cut $\IP^1_d$ along the  interval $I:=[-1,\infty]$ from $-1$ along the 
positive real axis.
Then the preimage of $I$ under $\pi$ will be a graph in $\IP^1_u$ with
vertices at each point of $\pi^{-1}(\{x,0,1,\infty\})$. The complement of the
graph will be the union of exactly $N$ connected components which are each
mapped isomorphically by $\pi$ onto $\IP^1_d\setminus I$, and in
particular the boundary of each component is the same as the boundary
of $\IP^1_d\setminus I$.
These connected components correspond to the branches of $\pi$ and the graph
specifies how to glue them together. In particular the graph determines the
permutation representation of the cover, since it shows us how to lift loops
in the base to paths in $\IP^1_u$; we just cross the corresponding edges
upstairs, and note which connected component we end up in.

Thus we need to draw the graphs in $\IP^1_u$. 
There are four types of vertices, depending on
if they lie over $-1,0,1,\infty$,
which we could draw as circles, squares, blobs and stars (say) respectively.
The number of each type of vertices is just the number of points of 
$\IP^1_u$ lying over the corresponding point amongst $-1,0,1,\infty$.
The corresponding ramification indices give the number of edges 
coming out of each vertex to each of 
the neighbouring vertices, and our task is to
join these edges together in a consistent manner.

For example for the first row of the above list, there should be $10$
branches and,
by examining the ramification indices, 
we see we need to draw a graph on $\IP^1_u$ out of the following pieces:

$\bullet$ $8$ circles with $1$ edge emanating from each, 
and $1$ circle with $2$ edges,

$\bullet$ $5$ squares with $4$ edges,

$\bullet$ $3$ blobs with $6$ edges and one blob with $2$ edges, and

$\bullet$ $1$ star with $7$ edges and $3$ stars with $1$ edge.

\noindent
The graph should divide the sphere into $10$ pieces and:

$\bullet$ Each edge from a circle should connect to a square,

$\bullet$ 
Two edges from each square should connect to a circle and the other 
two should connect to a blob (and, going around the square, the edges should
alternate between going to circles and blobs),

$\bullet$
Similarly half the edges from each blob should connect to squares and, again
alternating, the other half should connect to stars,

$\bullet$ Finally each edge from a star should connect to one of the blobs.

We leave the reader to draw such a graph 
(there are $15$ possibilities).\footnote{To count the possibilities, 
one may use Theorem 7.2.1 in Serre's book \cite{Serre-tgt} 
to count the number of such permutation representations, and then divide by 
conjugation action of the symmetric group, carefully computing the stabiliser.
To find all possibilities we draw some and then apply the natural action of
the pure three-string braid group to see if we get them all---here
all $15$ are braid equivalent.}
Given any such graph we can write down the monodromy of the pulled back
Fuchsian system on $\IP^1_u$ in terms of that of the hypergeometric system 
downstairs. 
Here the projective monodromy downstairs is a $(2,3,7)$ triangle group:
$$\Delta_{237}\cong\langle \ a,b,c\ \bigl\vert \ 
a^2=b^3=c^7=cba=1\ \rangle$$
which can be realised as a subgroup of $\PSL_2(\IC)$ in various ways (the
standard representation into $\PSL_2(\IR)$ plus its two Galois conjugates,
lying in $\PSU_2$).

Puncture $\IP^1_u$ at the four exceptional vertices (namely the $3$
stars with $1$ edge and the blob with $2$ edges) and choose generators 
$l_1,\cdots, l_4$ of the
fundamental group of this punctured sphere, with $l_4\circ \cdots \circ l_1$
contractible.
Then we can compute the image under $\pi$ of each $l_i$ in 
$\IP^1_d\setminus\{0,1,\infty\}$ and thereby write the monodromy of the
pulled back system as words in $a,b,c\in\Delta_{237}$.
With one such graph we obtained:
$$
L_1=caca^{-1}c^{-1},\quad 
L_2=c,\quad 
L_3=c^{-1}a^{-1}cac,\quad 
L_4=c^{-3}bc^3$$
where $L_i$ is the projective monodromy around $l_i$.
By construction $L_4\cdots L_1=1$ in $\Delta_{237}$.
Now by choosing an embedding of $\Delta_{237}$
in $\PSL_2(\IC)$
we get $L_i\in \PSL_2(\IC)$ and 
we can lift each $L_i$ to a matrix $M_i\in G$, (and possibly negate 
$M_4$ to ensure $M_4\cdots M_1=1$).
We obtain the representation 
$\rho$ with data
$$\theta=(2/7,2/7,2/7,1/3)\qquad\si=(1/3,1/3,1/7).$$

This completes our task of producing a representation in a 
finite $\F_2$ orbit.
Now we can apply our previous method to
construct the corresponding \PVI\  solution. Immediately, by computing the 
$\F_2$ orbit of the conjugacy class of $\rho$, we find the solution has genus
$1$ and $18$ branches, and that the parameters are not equivalent to those of
any known solution. 
Moreover it turns out that Jimbo's asymptotic formula may be applied to $17$
of the $18$ branches, and the asymptotics on the remaining branch may be
obtained by  the lemma in section 7 of \cite{icosa}.
Using this we can get the solution polynomial 
$F$ explicitly from the Puiseux expansions, and then look for a 
parameterisation of $F$. The result is:
\begin{gather}
\text{$2,3,7$ solution, genus one, $18$ branches,
$(\theta_1,\theta_2,\theta_3,\theta_4)=(2/7, 2/7, 2/7, 1/3)$:}\notag\\
y=
\frac{1}{2}-{\frac { \left( 3\,{s}^{8}-2\,{s}^{7}-4\,{s}^{6}-204\,{s}^{5}
-536\,{s}^{4}-1738\,{s}^{3}-5064\,{s}^{2}-4808\,s-3199 \right) u}
{4\, \left( {s}^{6}+196\,{s}^{3}+189\,{s}^{2}+756\,s+154 \right)  \left( {
s}^{2}+s+7 \right)  \left( s+1 \right) }}\notag \\
t=
\frac{1}{2}-{\frac { \left( {s}^{9}-84\,{s}^{6}-378\,{s}^{5}
-1512\,{s}^{4}-5208\,{s}^{3}-7236\,{s}^{2}-8127\,s-784 \right) u}
{432\,s \left( s+1 \right) ^{2} \left( {s}^{2}+s+7 \right) ^{2}}}
\label{eq: 237 soln}
\end{gather}
where 
$$u^2=s\,(s^2+s+7).$$

This solution is noteworthy in that currently there is no known relation to a
Fuchsian system with finite monodromy group (one might speculate as to the
existence of another realisation of \PVI\ in which this solution corresponds
to such a Fuchsian system, but this is unknown).

For the other three entries on the excerpt of
Doran's list above, we do not seem
to get new solutions, but it is interesting to identify them in any case.

The second entry, a family of degree $12$ covers, turns out to give the Klein
solution of \cite{k2p}.
The explicit family of covers has been found more recently 
in \cite{Kitaev-dessins} p.27.
There are $7$ different graphs one could draw, one of which is symmetrical and
they are all braid equivalent.
Using one of these graphs we obtain, as above, the  words:
$$
L_1=c^{-1}a^{-1}cac,\quad
L_2=c^3aca^{-1}c^{-3},\quad
L_3=c^2aca^{-1}c^{-2},\quad
L_4=a^{-1}c^{-1}a^{-1}c^2aca.$$
Choosing an appropriate embedding of $\Delta_{237}$ and lifting to 
$G$ we obtain the representation $\rho$ specified in \eqref{eq: klein data}.
In particular this gives a convenient way to prove that the projective
monodromy group of the family of $2\times 2$ 
Fuchsian systems determined by the Klein
solution is $\Delta_{237}$.
We just need to show that the $L_i$ generate all of $\Delta_{237}$, which we
will do in appendix \ref{apx: 237} below.

The third entry indicates a family of degree $12$ covers along which one
should pull back the $(2,3,8)$ triangle group.
This time there are $7$ graphs one could draw but they are not all braid
equivalent, there are two $P_3$ orbits, distinguished by the fact that
the monodromy group of the cover
is either $\Sym_{12}$ or a group of order $1536$.
For the degenerate case one finds the \PVI\ solution has just
two branches and is $y=t\pm\sqrt{t(t-1)}$ with parameters
$\theta=(1,1,1,7)/8$. (This is just the transform of the square root 
solution $y=\sqrt{t}$
under the Okamoto transformation $(y,t)\mapsto
(\frac{y-t}{1-t},\frac{t}{t-1})$).
The other case is more interesting; for one graph we obtain:
$$
L_1= aca^{-1},\quad
L_2=c^{-2}a^{-1}cac^{2},\quad
L_3=caca^{-1}c^{-1},\quad
L_4=a^{-1}caca^{-1}c^{-1}a$$
where now $a,b,c$ generate the $(2,3,8)$ triangle group:
$$\Delta_{238}\cong\langle \ a,b,c\ \bigl\vert \ 
a^2=b^3=c^8=cba=1\ \rangle.$$
Now we can choose an  embedding of $\Delta_{238}$ into $\PSL_2(\IC)$ and a
lift to $G$ (negating $M_4$ if necessary)
such that we obtain the representation $\rho$ 
with data 
$$\theta=(3/8, 3/8, 3/8, 5/8), \qquad\sigma=(1/2, 2/3, 1/3).$$
This is precisely that obtained above by applying an Okamoto transformation
to the $16$ branch octahedral solution (and gives a convenient way to prove,
in appendix \ref{apx: 237},
that the projective monodromy group is $\Delta_{238}$).

Finally there are $9$ graphs corresponding to the
last row of Doran's list, all braid equivalent.
Even though the graphs are the most complicated in this case (and there is a
quite attractive one with $4$-fold symmetry), this case  leads 
again 
to the $2$-branch \PVI\ solution $y=t\pm\sqrt{t(t-1)}$ with the parameters
$\theta=(1,1,1,6)/7$ (and $\si= (1/2,1/2,5/7)$ on one branch).

In conclusion we should mention that we do not know any other
finite $\F_2$ orbits of triples of elements of $\SL_2(\IC)$
(e.g. up to isomorphism as `sets with $\F_2$-action'); so far
they all either come from a finite subgroup or one of the two $(2,3,7)$ 
cases (the Klein solution or the genus one solution above).

\end{section}

\appendix
\begin{section}{} \label{apx: residues}
Here are the explicit formula from \cite{JM81} for the residue matrices $A_i$, 
of the isomonodromic family of Fuchsian systems corresponding to a
\PVI\ solution $y(t)$ with parameters $\theta_1,\ldots,\theta_4$.
The matrix entries are rational functions of 
$y,t,y'=\frac{dy}{dt},\{\theta_i\}$.
(Our coordinate $x$ is denoted $\wt z$ 
in \cite{JM81} and is related simply to $p$ which is the usual 
dual variable to $y=q$ in the Hamiltonian formulation  of \PVI. 
Also one should add 
$\diag(\theta_i,\theta_i)/2$ to our $A_i$ to obtain that of 
\cite{JM81}.)

$$A_i:=\left(\begin{matrix}
z_i+\theta_i/2 & -u_i z_i \\ 
(z_i+\theta_i)/u_i & -z_i-\theta_i/2 
\end{matrix}\right)\in \lsl_2(\IC)$$

where 
{\small
$$z_1=
y 
\frac{
E-k_2^2 (t+1)
}{t \theta_4},\qquad
z_2=
(y-t) 
\frac{
E+
t \theta_4 (y-1){x}
k_2^2-t k_1 k_2
}{t (t-1) \theta_4}
$$$$
z_3=
-(y-1) 
\frac{
E+\theta_4 (y-t) {x}
-k_2^2 t - k_1 k_2
}{(t-1) \theta_4},
$$

$$x=
p-\frac {\theta_1}{y}-\frac{\theta_2}{y-t}-\frac {\theta_3}{y-1},\quad
2 p=
\frac {\theta_1+\left( t-1 \right)y'}{y}+
\frac{\theta_2-1+y'}{y-t}+
\frac {\theta_3-t\,y'}{y-1}
$$

$$
u_1=\frac{y}{t z_1},\qquad
u_2=\frac{y-t}{t (t-1) z_2},\qquad
u_3=-\frac{y-1}{(t-1) z_3}$$
$$E=
y (y-1) (y-t) {x}^2+
\bigl(\theta_3 (y-t)+t \theta_2 (y-1)-2 k_2 (y-1) (y-t)\bigr) {x}
+k_2^2 y-k_2 (\theta_3+t \theta_2)$$
$$
k_1=(\theta_4-\theta_1-\theta_2-\theta_3)/2,\qquad
k_2=(-\theta_4-\theta_1-\theta_2-\theta_3)/2.$$
}

\end{section}
\begin{section}{} \label{apx: 237}
\begin{prop}
The $2\times 2$ Fuchsian systems corresponding to the Klein solution and to
the $18$ branch genus $1$ solution of section \ref{sn: img} have projective
monodromy group isomorphic to $\Delta_{237}$, and that corresponding to
the transformation of the $16$ branch octahedral appearing in 
section \ref{sn: img} has projective
monodromy group isomorphic to $\Delta_{238}$.
\end{prop}
{\noindent {\bf Proof.}\quad }
Since in section  \ref{sn: img} the projective 
monodromy groups were expressed as words in the generators of the respective
triangle groups, it is sufficient to check in each case that these words 
are in fact generators.
To do this we will repeatedly use the fact that in the group
$$\Delta=\langle \ a,b,c\ \bigl\vert\  a^2=b^3=c^n=cba=1\ \rangle$$
one has $c^b=bc^{-1}$ and $c^{b^{-1}}=c^{-1}b$, where in general we write
$x^y$ for $y^{-1}xy$. (These are easily verified,  for example the first is
true since $b^{-1}cbc=b^{-1}(cb)(cb)b^{-1}=b^{-2}=b$, using the fact that 
$a=a^{-1}=cb$.)
In particular we immediately deduce
$\Delta=\langle c,c^b\rangle=\langle c,c^{b^{-1}}\rangle$.

Now each case is an easy exercise. For the Klein case
we need to show $\langle L_i, i=1,2,3\rangle =\Delta$ where
$n=7$ and
$$
L_1=c^{ac},\quad
L_2=c^{ac^{-3}},\quad
L_3=c^{ac^{-2}}.$$
Up to conjugacy in $\Delta$, we have 
$\langle L_2,L_3\rangle = \langle p,c \rangle$
where $p=c^{ac^{-1}a}$. However, using $a=cb$ we see
$p=c^{cb^2}=c^{b^2}=c^{b^{-1}}$ so we are done.

For the other $(2,3,7)$ case corresponding to the genus one solution we have
$$L_1=c^{ac^{-1}},\quad 
L_2=c,\quad 
L_3=c^{ac}.$$
Thus $\langle L_2,L_3\rangle = \langle c,c^a\rangle=
\langle c,c^{b}\rangle$ since $a=cb$.

For the $(2,3,8)$  case we have
$L_1= c^a,
L_2=c^{ac^{2}},
L_3=c^{ac^{-1}}.$
Up to conjugacy $\langle L_1,L_3\rangle= \langle c, c^{ac^{-1}a}\rangle$
and as in the Klein case above $c^{ac^{-1}a}=bcb^{-1}$.\hfill$\square$
\end{section}

\begin{section}{} \label{apx: cdets}

At the request of the editors and of A. Kitaev, 
we will add some remarks to aid the reader interested in reproducing the
results of this article.
The main results are of two types:
1) classification of \PVI\  solutions coming from the binary tetrahedral and
octahedral groups and 2) construction of explicit \PVI\ solutions using
Jimbo's asymptotic formula.
For both 1) and 2) the details are parallel to those described  in
\cite{icosa} and \cite{k2p} resp., 
with the precise references as in the body of this article.
For 1) there are 3 steps:

\noindent$\bullet$
Prove that the relation of Okamoto equivalence is sandwiched between 
the relations of geometric and parametric equivalence, i.e. in symbols one has 
$\text{GE} \Rightarrow \text{OE} \Rightarrow \text{PE}. $
The second arrow is immediate by definition 
and the first is proved in Lemma 9 of \cite{icosa}.

\noindent$\bullet$ 
Compute the parameter equivalence classes in the set of parameters coming
from triples of generators of either the tetrahedral or octahedral group. 
This is as in section 3 of \cite{icosa}.
One first writes down the set of possible parameters $\theta$. 
This is a finite subset of $\IQ^4\subset\IR^4$.
Then one uses a
simple algorithm to move each of these points 
into a chosen affine $F_4$ 
alcove, using the standard action of the affine Weyl group $W_a(F_4)$ on
$\IR^4$
(this is entirely standard and the details are written in \cite{icosa}
Proposition 6). 
Then we count the number of distinct alcove points obtained.
By definition this is the number of ``parameter equivalence classes''.

\noindent$\bullet$ 
Compute the geometric equivalence classes in the set of linear representations
$\rho$ 
coming from either the binary 
tetrahedral or octahedral group. This amounts
to computing the orbits of an explicit action of a group on a finite set
(of size $520$, $3360$ resp.) and is carefully described in section 4 of
\cite{icosa}.

Some 
confidence that there is no computational error comes from the fact that
the geometric and parametric equivalence
classes turn out to coincide in both the cases considered here. 
Also Hall's formulae
\cite{Hall36} (computing
the number of generating triples) gives confidence that all the generating
triples have been computed correctly---since we get the right number 
of them. (In principle one can go through all triples of elements of the
finite group $\Gamma\subset \SL_2(\IC)$ and
throw out those that do not generate $\Gamma$. 
In the two cases at hand this is
feasible, but some simple tricks are useful in the icosahedral case.) 

Now we will move on to 2), constructing the solutions.
The main steps of the procedure used are as in \cite{k2p} 
(see especially p.193).
However with experience various tricks have been developed to speed up the
computation, so we will also 
detail some of these below (they are inessential if
one has a fast enough computer, as presumably future readers will have). 
The underlying strategy is analogous to that used in \cite{DubMaz00} although 
we do not in fact use any of their results.
(It was particularly troublesome to get the correct form of Jimbo's formula
in \cite{k2p},
which is the main ingredient and was not used in \cite{DubMaz00}.)

The basic steps are as follows:

1) We start with a linear representation $\rho$ living in a finite mapping
   class group orbit. The conjugacy class of $\rho$ is encoded in the
   seven-tuple 
$$\theta_1,\theta_2,\theta_3,\theta_4,\sigma_{12},\sigma_{23},\sigma_{13}.$$
Specifying these seven numbers is equivalent to specifying the numbers
$m_i=2\cos(\pi \theta_i),$ $m_{ij}=2\cos(\pi \si_{ij})$ provided we agree
 $\theta_i,\sigma_{ij}\in [0,1]$.
We compute the orbit of this $7$-tuple under 
the pure mapping class group $\cong 
\cF_2$. 
The formula for this action is given in \cite{icosa} section 4
(cf. also \eqref{eq: braiding} above).
This gives a list, of length $N$ say, 
of $7$-tuples, one for each branch of the corresponding \PVI\
solution. The values of the $\theta$'s will not vary on different branches
so the branches are parameterised by the values of the sigmas.
Let their values on the $k$th branch be denoted $\sigma^k_{ij}$,
$k=1,\ldots,N$.

2) Plug each $7$-tuple into Jimbo's asymptotic formula
(in the form in \cite{k2p} Theorem 4). 
This gives $N$ leading terms $y_k= a_k t^{b_k}+\cdots$ for $k=1,\ldots,N$
of the Puiseux expansion at $0$ 
of the \PVI\ solution $y(t)$ on the $N$ branches.
One will have $b_k=1-\sigma^k_{12}$ but $a_k$ is given by a 
very complicated, but explicit, formula.
(Jimbo's formula is not always applicable---cf. the discussion of `good'
solutions in \cite{icosa}, but often there is an equivalent solution for which
Jimbo's formula can be applied on every branch, or there is a degeneration of
Jimbo's formula (as in \cite{DubMaz00} 
or \cite{icosa} Lemma 19) which will compute
the remaining leading terms.)

3) Compute lots of terms in 
the Puiseux expansions of the solutions $y(t)$ on each branch.
These will be expansions in $s=t^{1/d_k}$ where $d_k$ is the denominator of
$b_k$ (when written in lowest terms).
Geometrically $d_k$ is the number of branches of $y$ that meet the given
branch  over $t=0$, 
i.e. it is the cycle length of the given branch in the permutation
representation of the solution curve as a cover of the $t$-line.
The expansions are computed recursively by substituting back into the \PVI\ 
equation; at each step this leads to a linear equation for the coefficient of 
the next term in the expansion.

4) Use these expansions to determine the coefficients of the solution 
polynomial $F(y,t)$. (This determines $y$ as an algebraic function of $t$
by the condition $F(y,t)=0$.)
Since $F$ is a  polynomial (of degree $N$ in $y$) 
this is clearly possible since we have arbitrarily many 
terms of each Puiseux expansion; $F(y_k(s),s^{d_k})=0$ 
for all $s$ and for each branch $y_k$ of the solution.
(Thus 
in principle just one expansion is needed, 
not the expansion for all branches.)
Given $F(y,t)$ one may check symbolically that it specifies a solution to
\PVI, using implicit differentiation.

5) Compute a parameterisation of the resulting curve $F(y,t)=0$.
(As mentioned in the acknowledgments the author is grateful to 
Mark van Hoeij for help with this last step.)
In general this will be simpler to write down than the polynomial $F$.

\ 

Now we will list some of the tricks we have found useful in carrying out the
above steps.

1) One needs to convert the numbers $a_k$ given by Jimbo's formula into
   algebraic numbers. In the examples so far 
   this can be done by raising $a_k$ (and/or
   its real/imaginary parts) to the $d_k$-th power until a rational number is
   obtained (which can be ascertained by looking at continued fraction
   expansions). 

2) Reduce the number of Puiseux expansions to compute:
the $d_k$ branches which meet the given branch over zero will have Galois
conjugate expansions. These can be obtained from one another by multiplying
$s$ by a $d_k$-th root of unity. Also, when choosing which of these $d_k$
expansions to actually compute it is good to choose the real one, if possible.
(Also sometimes some expansions are complex conjugate to others so further
optimisations are possible.)

3) Reduce the degree of the field extension used to compute the expansion:
In computing the Puiseux expansions one is often working over a finite
extension of $\IQ$, such as $\IQ(6^{1/7})$.
Often the degree of this extension can be reduced by taking the expansion in a
variable $h=c\times s$ for a suitable constant $c$, rather than in 
$s=t^{1/d_k}$.  This trick was very useful for computing the larger solutions
(with $\ge 15$ branches say).

4) To obtain the coefficients of the polynomial $F$ from the Puiseux expansion
   we use the trick suggested in \cite{DubMaz00}:
Write $F$ in the form
$$F=q(t)y^N+p_{N-1}(t)y^{N-1}+\cdots+p_1(t)y+p_0(t)$$
for polynomials $p_i,q$ in $t$ and
define rational functions $r_i(t):=p_i/q$ for $i=0,\ldots,N-1$.
If $y_1,\ldots,y_N$ denote the (locally defined) solutions on the branches
then for each $t$ we have that $y_1(t),\ldots,y_N(t)$ are the roots of 
$F(t,y)=0$ and it follows that
$$y^N+r_{N-1}(t)y^{N-1}+
\cdots+r_1(t)y+r_0(t)=(y-y_1(t))(y-y_2(t))\cdots(y-y_N(t)).$$
Thus, expanding the product on the right, the rational functions $r_i$ are
obtained as symmetric polynomials in the $y_i$:
$$r_0=(-1)^N y_1\cdots y_N,\qquad \cdots\qquad ,r_{N-1}=-(y_1+\cdots+y_N).$$ 
Since the $r_i$ are global rational functions, the Puiseux expansions of the
$y_i$ give the Laurent expansions at $0$ of the $r_i$.
Clearly only a finite number of terms of each Laurent expansion 
are required to determine each $r_i$, and indeed it is simple to convert these
truncated Laurent expansions into global rational functions.
(This is easily done by Pad\'e approximation, e.g. as implemented in the Maple
command ``$\text{convert}(\ ,\text{ratpoly})$''.)

5) Much time may be saved by carefully choosing the 
   representative for the solution in the first place (i.e. try to 
   choose an equivalent solution for which the polynomial $F$ is simpler).
   Heuristically this can be estimated by seeing how complicated the algebraic
   numbers $a_k$ are (or by seeing how complicated the coefficients of the
   polynomial $q(t)$ are; this is usually 
   easily obtained from $(y_1+\cdots+y_N)$,
   i.e. before having to compute complicated symmetric functions).

6) Use Okamoto symmetries wherever possible: e.g. if (we can arrange that)
   the solution has the
   symmetry $(y,t)\mapsto (y/t,1/t)$, swapping $\theta_2$ and $\theta_3$
then the coefficients of each $p_i,q$ should be symmetrical, thereby
   essentially halving
   the number of coefficients that need to be computed.
(Also for the $24$-branch icosahedral solution in \cite{icosa} it was too
   cumbersome to compute the longest symmetric functions, corresponding to the
   `middle' polynomials $p_i$, but, by using another Okamoto symmetry, the
   outstanding coefficients could be determined in terms of those we were able
   to compute.)

7) Finally there are various optimisations that can be made (especially in
   computing the symmetric functions of the Puiseux expansions) if we expect
   $F$ to have integer coefficients 
(which is the case for all examples so far).

\ 

\renewcommand{\baselinestretch}{1}           
\Small
{\bf Acknowledgments.}\ \ 
The author is very grateful to 
Mark van Hoeij for help computing the more difficult
parameterisations of the curves $F=0$, 
and to both C. Doran and A. Kitaev for explaining various aspects of their
work. 

\ 

After this work was complete A. Kitaev informed the author that he had found
an explicit family of covers corresponding to the genus one $(2,3,7)$
solution of section \ref{sn: img} and had obtained a similar solution.
Happily the solution here and that of Kitaev are not related by 
Okamoto transformations, but arise by choosing different embeddings of
$\Delta_{237}$ into $\PSL_2(\IC)$. In fact there are three inequivalent 
choices, 
corresponding to the three conjugacy classes of order $7$ elements in
$\PSL_2(\IC)\cong \SO_3(\IC)$.
(This is analogous 
to the sibling solutions of \cite{icosa} which arose from the two
classes of order $5$ elements.)
The third inequivalent \PVI\ 
solution is:
{\tiny
$$y= \frac{1}{2}-
{\frac { \left( {s}^{10}+5\,{s}^{9}+24\,{s}^{8}+20\,{s}^{7}-
266\,{s}^{6}-2874\,{s}^{5}-14812\,{s}^{4}-40316\,{s}^{3}-85359\,{s}^{2
}-100067\,s-67396 \right) u}{ 16\left( s+1 \right)  \left( {s}^{2}+s+7
 \right)  \left( 5\,{s}^{6}+63\,{s}^{5}+252\,{s}^{4}+854\,{s}^{3}+1449
\,{s}^{2}+1827\,s+2030 \right) }}
$$
}
with $t,u,s$ exactly as in \eqref{eq: 237 soln} and 
$\theta=(4/7,4/7,4/7,1/3)$.

\end{section}

\renewcommand{\baselinestretch}{1}              %
\normalsize
\bibliographystyle{amsplain}    \label{biby}
\bibliography{../thesis/syr}    

\providecommand{\bysame}{\leavevmode\hbox to3em{\hrulefill}\thinspace}
\providecommand{\MR}{\relax\ifhmode\unskip\space\fi MR }
\providecommand{\MRhref}[2]{%
  \href{http://www.ams.org/mathscinet-getitem?mr=#1}{#2}
}
\providecommand{\href}[2]{#2}
\begin{thebibliography}{10}

\bibitem{And-Kit-CMP}
F.~V. Andreev and A.~V. Kitaev, \emph{Transformations {$RS\sp 2\sb 4(3)$} of
  the ranks {$\leq 4$} and algebraic solutions of the sixth {P}ainlev\'e
  equation}, Comm. Math. Phys. \textbf{228} (2002), no.~1, 151--176.

\bibitem{AnoBol94}
D.~V. Anosov and A.~A. Bolibruch, \emph{The {R}iemann-{H}ilbert problem},
  Aspects of Mathematics, E22, Friedr. Vieweg \& Sohn, Braunschweig, 1994.

\bibitem{Bald-Dwork-algsols}
F.~Baldassarri and B.~Dwork, \emph{On second order linear differential
  equations with algebraic solutions}, Amer. J. Math. \textbf{101} (1979),
  no.~1, 42--76.

\bibitem{B-vdW}
F.~Beukers and A.~van~der Waall, \emph{Lam\'e equations with algebraic
  solutions}, J. Differential Equations \textbf{197} (2004), no.~1, 1--25.

\bibitem{icosa}
P.~P. Boalch, \emph{The fifty-two icosahedral solutions to {P}ainlev\'e {VI}},
  J. Reine Angew. Math., to appear, (math.AG/0406281, v.7).

\bibitem{ipc}
\bysame, \emph{{Higher genus icosahedral {P}ainlev\'e curves}},
  math.AG/0506407.

\bibitem{k2p}
\bysame, \emph{From {K}lein to {P}ainlev\'e via {F}ourier, {L}aplace and
  {J}imbo}, Proc. London Math. Soc. \textbf{90} (2005), no.~3, 167--208.

\bibitem{Chuck1}
C.~F. Doran, \emph{Algebraic and geometric isomonodromic deformations}, J.
  Differential Geom. \textbf{59} (2001), no.~1, 33--85.

\bibitem{Dub95long}
B.~Dubrovin, \emph{Geometry of 2{D} topological field theories}, Integrable
  Systems and Quantum Groups (M.Francaviglia and S.Greco, eds.), vol. 1620,
  Springer Lect. Notes Math., 1995, pp.~120--348.

\bibitem{DubMaz00}
B.~Dubrovin and M.~Mazzocco, \emph{Monodromy of certain {P}ainlev\'e-{V}{I}
  transcendents and reflection groups}, Invent. Math. \textbf{141} (2000),
  no.~1, 55--147. \MR{2001j:34114}

\bibitem{Hall36}
P.~Hall, \emph{The {E}ulerian functions of a group}, Quart. J. Math. Oxford
  Ser. 7 (1936), 134--151.

\bibitem{Hit-Poncelet}
N.~J. Hitchin, \emph{Poncelet polygons and the {P}ainlev\'e equations},
  Geometry and analysis (Bombay, 1992), Tata Inst. Fund. Res., Bombay, 1995,
  pp.~151--185. \MR{97d:32042}

\bibitem{Hit-Octa}
\bysame, \emph{A lecture on the octahedron}, Bull. London Math. Soc.
  \textbf{35} (2003), 577--600.

\bibitem{IIS}
M.~Inaba, K.~Iwasaki, and M.-H. Saito, \emph{B\"acklund transformations of the
  sixth {P}ainlev\'e equation in terms of {R}iemann-{H}ilbert correspondence},
  I. M. R. N. (2004), no.~1, 1--30, math.AG/0309341.

\bibitem{Jimbo82}
M.~Jimbo, \emph{Monodromy problem and the boundary condition for some
  {P}ainlev\'e equations}, Publ. Res. Inst. Math. Sci. \textbf{18} (1982),
  no.~3, 1137--1161. \MR{85c:58050}

\bibitem{JM81}
M.~Jimbo and T.~Miwa, \emph{Monodromy preserving deformations of linear
  differential equations with rational coefficients {II}}, Physica 2D (1981),
  407--448.

\bibitem{Katz-Algsols}
N.~M. Katz, \emph{Algebraic solutions of differential equations
  ({$p$}-curvature and the {H}odge filtration)}, Invent. Math. \textbf{18}
  (1972), 1--118.

\bibitem{Kitaev-dessins}
A.~V. Kitaev, \emph{Dessins d'enfants, their deformations and algebraic the
  sixth {P}ainlev\'e and {G}auss hypergeometric functions}, nlin.SI/0309078,
  v.3.

\bibitem{Kitaev-quad-p6}
\bysame, \emph{Quadratic transformations for the sixth {P}ainlev\'e equation},
  Lett. Math. Phys. \textbf{21} (1991), no.~2, 105--111.

\bibitem{Kit-sfit6}
\bysame, \emph{Special functions of isomonodromy type, rational transformations
  of the spectral parameter, and algebraic solutions of the sixth {P}ainlev\'e
  equation}, Algebra i Analiz \textbf{14} (2002), no.~3, 121--139.

\bibitem{Klein-ico}
Felix Klein, \emph{Lectures on the icosahedron and the solution of equations of
  the fifth degree}, Dover Publications Inc., New York, N.Y., 1956.

\bibitem{Mal-imd1long}
B.~Malgrange, \emph{Sur les deformations isomonodromiques. {I}. singularites
  regulieres}, S\'eminaire E.N.S. Math\'ematique et Physique (Boston)
  (L.~Boutet~de Monvel, A.~Douady, and J.-L. Verdier, eds.), Progress in Math.,
  vol.~37, Birkh\"auser, 1983, pp.~401--426.

\bibitem{SMJ}
M.~Sato, T.~Miwa, and M.~Jimbo, \emph{Holonomic quantum fields {I-V}}, Publ.
  RIMS Kyoto Univ. \textbf{14,15,15,15,16} (1978,1979,1979,1979,1980),
  223--267,201--278,577--629,871--972,531--584, resp.

\bibitem{Schles}
L.~Schlesinger, \emph{\"{U}ber eine {K}lasse von {D}ifferentialsystemen
  beliebiger {O}rdnung mit feten kritischen {P}unkten}, J. f\"ur Math.
  \textbf{141} (1912), 96--145.

\bibitem{Schwarz}
H.~A. Schwarz, \emph{{\"U}ber diejenigen {F}\"alle, in welchen die {G}aussische
  hypergeometrische reihe eine algebraische {F}unktion ihres vierten {E}lements
  darstellt}, J. Reine Angew. Math. \textbf{75} (1873), 292--335.

\bibitem{Serre-tgt}
J.-P. Serre, \emph{Topics in {G}alois theory}, Research Notes in Mathematics,
  vol.~1, Jones and Bartlett Publishers, Boston, MA, 1992.

\bibitem{vdP-Ulmer-defg}
M.~van~der Put and F.~Ulmer, \emph{Differential equations and finite groups},
  J. Algebra \textbf{226} (2000), no.~2, 920--966.

\bibitem{watanabePVI}
H.~Watanabe, \emph{Birational canonical transformations and classical solutions
  of the sixth {P}ainlev\'e equation}, Ann. Scuola Norm. Sup. Pisa Cl. Sci. (4)
  \textbf{27} (1998), no.~3-4, 379--425.

\end{thebibliography}
\end{document}